\documentclass[12pt]{article}

\usepackage[dvipdfm, bookmarksopen=true]{hyperref}
\usepackage[centertags]{amsmath}
\usepackage{amsfonts}
\usepackage{amssymb}
\usepackage{amsthm}
\usepackage{bm}

\theoremstyle{plain}
\newtheorem{thm}{Theorem}[section]
\newtheorem{cor}[thm]{Corollary}
\newtheorem{lem}[thm]{Lemma}
\newtheorem{prop}[thm]{Proposition}

\theoremstyle{definition}
\newtheorem{defn}[thm]{Definition}
\newtheorem{exam}[thm]{Example}
\newtheorem{rem}[thm]{Remark}
\theoremstyle{remark}
\numberwithin{equation}{section}

\newcommand{\beast}{\begin{eqnarray*}}
\newcommand{\eeast}{\end{eqnarray*}}

\title{The edge-flipping group of a graph\thanks{Research partially supported by the NSC grant 97-2115-M-009-002 of
Taiwan R.O.C.}}

\date{June 12, 2009}

\author{Hau-wen Huang\footnote{Corresponding author}\and Chih-wen Weng\footnote{$~^\dag$Department of Applied Mathematics National Chiao Tung University 1001 Ta Hsueh
Road Hsinchu, Taiwan 300, R.O.C.}}

\begin{document}
\maketitle

\begin{abstract}
Let $X=(V,E)$ be a finite simple connected graph with $n$ vertices and $m$ edges. A configuration is an assignment of
one of two colors, black or white, to each edge of $X.$ A move applied to a configuration is to select a black edge $\epsilon\in E$ and change the colors of all adjacent edges of $\epsilon.$ Given an initial configuration and a final configuration, try to find a sequence of moves that transforms the initial configuration into the final configuration. This is the edge-flipping puzzle on $X,$ and it corresponds to a group action. This group is called the edge-flipping group $\mathbf{W}_E(X)$ of $X.$ This paper shows that if $X$ has at least three vertices, $\mathbf{W}_E(X)$ is isomorphic to a semidirect product of $(\mathbb{Z}/2\mathbb{Z})^k$ and the symmetric group $S_n$ of degree $n,$ where
$k=(n-1)(m-n+1)$ if $n$ is odd, $k=(n-2)(m-n+1)$ if $n$ is even, and $\mathbb{Z}$ is the additive
group of integers.
\end{abstract}

{\footnotesize{\it Keywords:} group actions; orbits; semidirect products.}

\section{Introduction}\label{Introduction}
An ordered pair $X=(V,E)$ is a {\it finite simple graph} if $V$ is a finite
set and $E$ is a set of some 2-element subsets of $V.$ The
elements of $V$ are called {\it vertices} of $X$ and the elements of $E$ are called {\it edges} of $X$. Two vertices $u,v$ of $X$ are {\it neighbors} if $\{u,v\}\in E.$ A finite simple graph $X=(V,E)$ is {\it
connected} if for any two distinct vertices $u,v\in V$ there exists a subset
$\{\{u_0,u_1\},\{u_1,u_2\},\ldots,\{u_{k-1},u_k\}\}$ of $E$ with $u_0=u$ and $u_k=v.$

Throughout this paper, $X=(V,E)$ is a finite simple connected graph with $|V|=n$ and $|E|=m.$ This paper focuses on  two flipping puzzles defined on the graph $X$ as follows. A {\it configuration} of
the first puzzle (second puzzle, respectively) is an assignment of one of two colors, black or white, to each edge of $X$ (vertex of $X,$ respectively). A move applied to a configuration is to select a black edge $\epsilon\in E$ (black vertex $v\in V,$ respectively) and change the colors of all edges $\epsilon'$ with $|\epsilon'\cap\epsilon|=1$ (all neighbors of $v,$ respectively). Given an initial configuration and a final configuration, try to find a sequence of moves that transforms the initial configuration into the final configuration. This is the {\it edge-flipping puzzle} (the {\it vertex-flipping puzzle}, respectively) on $X.$ A set $O$ of some configurations is an {\it orbit} of the edge-flipping puzzle (the vertex-flipping puzzle, respectively) on $X$ if for any two configurations in $O,$ one can reach the other by a sequence of moves. The edge-flipping puzzle corresponds to a group action, and this group is called the edge-flipping group $\mathbf{W}_E(X)$ of $X=(V,E).$ See Section 3 for details.

The orbits of the edge-flipping puzzle on $X$ have been determined. If $X$ is a tree with at least three vertices, the edge-flipping group $\mathbf{W}_E(X)$ of $X$ is isomorphic to the symmetric group $S_n$ of degree $n.$ Wu \cite{w:08} illustrated both of these results. The main goal of the current study is to produce the complement part of his second result, and show that if $X$ has at least three vertices, then $\mathbf{W}_E(X)$ is isomorphic to
$$
\left\{
\begin{array}{ll}
(\mathbb{Z}/2\mathbb{Z})^{(n-1)(m-n+1)}\rtimes S_n, &\hbox{if $n$ is odd;}\\
(\mathbb{Z}/2\mathbb{Z})^{(n-2)(m-n+1)}\rtimes S_n, &\hbox{if $n$ is even,}
\end{array}
\right.
$$
where $\mathbb{Z}$ is the additive group of integers. Section 5 explains the structure of $\mathbf{W}_E(X)$ in greater detail. 

The development and history of vertex-flipping puzzles can be found in the literature \cite{aeb,mkc:04,mkc:06,hErik:01,hw:pre2,jww:08,xw:07,wu:06}. For example, the vertex-flipping puzzles implicitly appear in Chuah and Hu's papers \cite{mkc:04,mkc:06} when they study the equivalence classes of Vogan diagrams and extended Vogan diagrams \cite{pb:00,pb:02,k:96}. Note that the vertex-flipping puzzles are called lit-only $\sigma$-games in \cite{jww:08,xw:07,wu:06}.

The vertex-flipping puzzle on $X$ also corresponds to a group action \cite{xw:07}, and this group is called the vertex-flipping group $\mathbf{W}_{V}(X)$ of $X.$ Some properties of the vertex-flipping group $\mathbf{W}_{V}(X)$ of $X$ have been known. For example, $\mathbf{W}_{V}(X)$ has the trivial center, and $\mathbf{W}_{V}(X)$ is a homomorphic image of the Coxeter group $W$ of $X.$ If $X$ is a simply-laced Dynkin diagram, $\mathbf{W}_{V}(X)$ is isomorphic to the quotient group of $W$ by its center $Z(W);$ moreover, $|Z(W)|=1$ or $2.$ See \cite{hw:pre} for details. In other words, $\mathbf{W}_V(X)$ can be treated as a combinatorial version of the reflection groups on real vector spaces. Although J. Humphreys gave a faithful geometric representation of any Coxeter group $W$ in a real vector space \cite[Section 5.3]{h:90}, the group structures of $W$ and $\mathbf{W}_V(X)$ are worthy of further study.

On the other hand, the edge-flipping puzzle on $X$ is the vertex-flipping puzzle played on the line graph $L(X)$ of $X.$ The edge-flipping group $\mathbf{W}_E(X)$ of $X$ is also the vertex-flipping group of $L(X).$ The main result of this study can be used to find the group structures of the vertex-flipping groups of some graphs, which do not need to be line graphs. See Section 6 for details. Note that line graphs are classified in terms of nine forbidden induced subgraphs \cite{lb:68,west:01}. 

\section{Edge spaces and bond spaces}\label{edge_bond_spaces}
Let $X=(V,E)$ denote a finite simple connected graph with $|V|=n$ and $|E|=m.$ In this section we give some basic definitions and properties about the edge space and the bond space of $X.$ The reader may refer to \cite[p.23--p.28]{rd:05} for details. 
Let $\mathcal{E}$ denote the power set of $E.$ Let $\mathbf{F}_2=\{0,1\}$ denote the 2-element field. For $F,F'\in\mathcal{E}$, define $F+F':=\{\epsilon\in E~|~\epsilon\in F\cup F',\epsilon\notin F\cap F'\}$; i.e., the symmetric difference of $F$ and $F',$ and define $1\cdot F:=F$ and $0\cdot F:=\emptyset,$ the empty set. Then $\mathcal{E}$ forms a vector space over $\mathbf{F}_2$ and is called the {\it edge space} of $X.$ Note that the zero element of  $\mathcal{E}$ is $\emptyset$ and $-F=F$ for $F\in \mathcal{E}.$ Since $\{\{\epsilon\}~|~\epsilon\in E\}$ is a basis of $\mathcal{E},$ we have ${\rm dim}~\mathcal{E}=m.$ 
In the same way as above, the power set $\mathcal{V}$ of $V$ also forms a vector space over $\mathbf{F}_2$ with symmetric difference as vector addition, and we call $\mathcal{V}$ the {\it vertex space} of $X.$ Clearly, ${\rm dim}~\mathcal{V}=n.$

For a subset $U$ of $V,$ let $E(U)$ denote the subset of $E$ consisting of all edges of $X$ that have exactly one element in $U.$ In graph theory, $E(U)$ is often called an {\it edge cut} of $X$ if $U$ is a nonempty and proper subset of $V$. Let $E(v):=E(\{v\})$ for $v\in V$ and notice that $E(\epsilon)=E(\{x,y\})$ for $\epsilon=\{x,y\}\in E.$

\begin{prop}\label{bond_1}
Let $X=(V,E)$ be a finite simple connected graph. Then the following (i),(ii) hold.
\begin{enumerate}
\item[(i)]
Each $\epsilon=\{x,y\}\in E$ lies in exactly two edge cuts $E(x)$ and $E(y)$ among $E(v)$ for all $v\in V.$
\item[(ii)]
For $U\subseteq V,$ $E(U)=\sum_{v\in U}E(v).$
\end{enumerate}
\end{prop}
\begin{proof}
(i) is immediate from the definition of $E(v)$ for $v\in V.$ (ii) is immediate from (i) and the definition of $E(U).$
\end{proof}

The {\it bond space} $\mathcal{B}$ of $X$ is the subspace of $\mathcal{E}$ spanned by $E(v)$ for all $v\in V.$ In the following, we give some basic properties of $\mathcal{B}.$

\begin{prop}\label{dimB} Let $X=(V,E)$ be a finite simple connected graph with $n$ vertices, and let $\mathcal{B}$ be the bond space of $X.$ Then the following (i)-(iv) hold.
\begin{enumerate}
\item[(i)] $\mathcal{B}=\{E(U)~|~U\subseteq V\}.$

\item[(ii)] ${\rm dim}~\mathcal{B}=n-1.$

\item[(iii)] For $u\in V,$ $E(u)=\sum_{v\in V-\{u\}}E(v).$

\item[(iv)] For $u\in V,$ the set $\{E(v)~|~v\in V-\{u\}\}$ is a basis of $\mathcal{B}.$
\end{enumerate}
\end{prop}
\begin{proof}
(i) follows immediately from Proposition~\ref{bond_1}(ii). Note that the map from the vertex space $\mathcal{V}$ onto the bond space $\mathcal{B}$ of $X,$ defined by
$$
U\mapsto E(U){\rm~for~}U\in\mathcal{V},
$$
is a linear transformation with kernel $\{\emptyset, V\}.$ Hence ${\rm dim}~\mathcal{B}=n-1$ and this prove (ii). Let $u\in V.$  Since $E(V)=\emptyset,$ we have
\begin{align*}
E(u)&=E(u)+E(V)\\
&=\sum_{v\in V-\{u\}}E(v).
\end{align*}
This proves (iii). Also, the set $\{E(v)~|~v\in V-\{u\}\}$ is a basis of $\mathcal{B}$ since it has at most $n-1$ elements and spans $\mathcal{B}$ by (iii). This proves (iv).
\end{proof}

It is not hard to see that there exists an $(n-1)$-element subset $T$ of $E$ such that $(V,T)$ is connected since $X=(V,E)$ is connected. We call such $T$ a {\it spanning tree} of $E.$ The following proposition says that $\{F~|~F\subseteq E-T\}$ is a set of coset representatives of $\mathcal{B}$ in $\mathcal{E}.$

\begin{prop}\label{coset_rep}
Let $X=(V,E)$ be a finite simple connected graph with $n$ vertices and $m$ edges. Let $\mathcal{E}$ and $\mathcal{B}$ be the edge space and bond space of $X$ respectively. Then the subset $\{F~|~F\subseteq E-T\}$ of $\mathcal{E}$ is a set of coset representatives of $\mathcal{B}$ in $\mathcal{E},$ where $T$ is a spanning tree of $E.$
\end{prop}
\begin{proof}
Note that there are $2^{m-n+1}$ cosets of $\mathcal{B}$ in $\mathcal{E}$ because of ${\rm dim}~\mathcal{B}=n-1$ and ${\rm dim}~\mathcal{E}=m.$ It is clear that $|\{F~|~F\subseteq E-T\}|=2^{m-n+1}.$ For any two distinct $F,F'\subseteq E-T,$ the graph $(V,E-(F-F'))$ is still connected since $T\subseteq E-(F-F'),$ which implies that $F-F'$ is not an edge cut of $X;$ i.e., $F-F'\notin \mathcal{B}.$ Hence $\{F~|~F\subseteq E-T\}$ is a set of coset representatives of $\mathcal{B}$ in $\mathcal{E}.$
\end{proof}

\section{Edge-flipping groups and invariant subsets}
In this and next sections, we rephrase some results of \cite{w:08} in order to facilitate our work. Let $X=(V,E)$ denote a finite simple connected graph with $|V|=$ $n$ and $|E|=m.$ Let $\mathcal{E}$ and $\mathcal{B}$ denote the edge space and bond space of $X$ respectively. We regard every configuration of the edge-flipping puzzle on $X$ as an element $G$ of $\mathcal{E},$ where $G$ consists of all black edges. We give a new interpretation of the moves in the edge-flipping puzzle on $X$: on each round, we select an edge $\epsilon\in E.$ If $\epsilon$ is a black edge, then we change the colors of all edges $\epsilon'$ with $|\epsilon'\cap \epsilon|=1;$ otherwise, we do nothing. Clearly, the orbits of the edge-flipping puzzle of $X$ under the new moves are unchanged. However, the new move by selecting an edge $\epsilon$ of $X$ is corresponding to the map
$\bm{\rho}_{\epsilon}:\mathcal{E}\rightarrow\mathcal{E}$ defined by
\begin{align}\label{def_e}
\bm{\rho}_{\epsilon}G=\left\{
\begin{array}{ll}
G+E(\epsilon), & \hbox{if $\epsilon\in G;$}\\
G,     & \hbox{otherwise}
\end{array}
\right.
\end{align}
for $G\in \mathcal{E}.$ From the definition of $\bm{\rho}_{\epsilon},$ we see that $\bm{\rho}_{\epsilon}^2$ is the identity map on $\mathcal{E}.$ In particular,  $\bm{\rho}_{\epsilon}$ is invertible. It is straightforward to check that $\bm{\rho}_{\epsilon}$ is a linear transformation on $\mathcal{E}.$ Hence $\bm{\rho}_{\epsilon}$ is an element in the general linear group ${\rm GL}(\mathcal{E})$ of $\mathcal{E}.$

\begin{defn}\label{edge_flipping_gp}
The {\it edge-flipping group} $\mathbf{W}_E(X)$ of $X=(V,E)$ is the subgroup of the general linear group ${\rm GL}(\mathcal{E})$ of $\mathcal{E}$ generated by the set $\{\bm{\rho}_{\epsilon}~|~\epsilon\in E\}.$ 
\end{defn}

Recall that $I$ is an {\it invariant subset} of $\mathcal{E}$ under $\mathbf{W}_E(X)$ if $I\subseteq \mathcal{E}$ and $\mathbf{W}_E(X)I\subseteq I.$ In the following, we give some significant invariant subsets of $\mathcal{E}$ under $\mathbf{W}_E(X).$ 

\begin{prop}\label{coset}(\cite{w:08})
Let $X=(V,E)$ be a finite simple connected graph. Let $\mathcal{E}$ and $\mathcal{B}$ be the edge space and bond space of $X$ respectively. Then each coset of $\mathcal{B}$ in $\mathcal{E}$ is an invariant subset of $\mathcal{E}$ under $\mathbf{W}_E(X).$
\end{prop}
\begin{proof}
It suffices to show that $\bm{\rho}_{\epsilon}G$ is in $G+\mathcal{B}$ for any $\epsilon\in E$ and $G\in\mathcal{E}.$ By (\ref{def_e}), $\bm{\rho}_{\epsilon}G$ is equal to either $G+E(\epsilon)$ or $G.$ Hence, $\bm{\rho}_{\epsilon}G\in G+\mathcal{B}$ since $E(\epsilon)\in \mathcal{B}.$
\end{proof}

From now to the end of this section, we shall study the group action of $\mathbf{W}_E(X)$ on the bond space $\mathcal{B}$ of $X.$ Recall that the set $\{E(v)~|~v\in V\}$ spans $\mathcal{B}.$ We first determine the cardinality of the set $\{E(v)~|~v\in V\}$ as follows. 

\begin{lem}\label{deg_sym}
Let $X=(V,E)$ be a finite simple connected graph with $|V|=n.$  Then
$$
|\{E(v)~|~v\in V\}|=\left\{
\begin{array}{ll}
n, &\hbox{if $n\geq 3;$}\\
1, &\hbox{otherwise.}
\end{array}
\right.
$$
\end{lem}
\begin{proof}
If $X$ is a one-vertex graph, then $\{E(v)~|~v\in V\}=\{\emptyset\},$ and if $X$ is a connected graph with two vertices, say $x$ and $y,$ then $E(x)=E(y)$ and hence $|\{E(v)~|~v\in V\}|=1.$ Now suppose $n\geq 3.$ Pick two distinct vertices $u,v\in V,$ we show $E(u)\not=E(v).$ Since the edge cut $E(\{u,v\})$ is nonempty and by Proposition~\ref{bond_1}(ii), $E(u)+E(v)=E(\{u,v\})\not=\emptyset;$ i.e., $E(u)\not=E(v).$
\end{proof}

Throughout the remainder of this section, we assume $n\geq 3$ and we denote the symmetric group on $\{E(v)~|~v\in V\}$ of degree $n$ by $S_n.$ Suppose $\epsilon=\{x,y\}\in E.$  From Proposition~\ref{bond_1}(i) and the definition of $\bm{\rho}_{\epsilon}$ in (\ref{def_e}), we know that $\bm{\rho}_{\epsilon}$ fixes all $E(v)$'s except $E(x)$ and $E(y).$ Also from Proposition~\ref{bond_1}(ii), we know that $E(\epsilon)=E(x)+E(y),$ and hence $\bm{\rho}_{\epsilon}E(x)=E(y)$ and $\bm{\rho}_{\epsilon}E(y)=E(x).$ In brief, the mapping of $\bm{\rho}_{\epsilon}$ on $\{E(v)~|~v\in V\}$ is the transposition $(E(x),E(y))$ in $S_n.$ Since each element $\mathbf{g}\in\mathbf{W}_E(X)$ is generated by $\bm{\rho}_{\epsilon}$ for $\epsilon\in E,$ the mapping of $\mathbf{g}$ on $\{E(v)~|~v\in V\}$ is like a permutation in $S_n.$ Hence we have the following definition.

\begin{defn}\label{alpha}
Let $\alpha:\mathbf{W}_E(X)\rightarrow S_n$ denote the group homomorphism defined by
$$
\alpha(\mathbf{g})(E(v))=\mathbf{g}E(v)
$$
for $v\in V$ and $\mathbf{g}\in \mathbf{W}_E(X).$
\end{defn}

Let $T$ be a spanning tree of $E,$ and let $\mathbf{W}_E(X)_T$ denote the subgroup of $\mathbf{W}_E(X)$ generated by the set $\{\bm{\rho}_{\epsilon}~|~\epsilon\in T\}.$ We say that $X$ is a {\it tree} if $E=T.$ The following lemma shows that $\mathbf{W}_E(X)$ is isomorphic to $S_n$ if $X$ is a tree.

\begin{lem}(\cite[Theorem 8]{w:08})\label{WWTontoSn}
Let $X=(V,E)$ be a finite simple connected graph with $|V|=n\geq 3.$ Let $S_n$ be the symmetric group on $\{E(v)~|~v\in V\}$ of degree $n.$ Then $\alpha(\mathbf{W}_E(X))=\alpha(\mathbf{W}_E(X)_T)=S_n,$ where $T$ is a spanning tree of $E.$  Moreover,
$\mathbf{W}_E(X)$ is isomorphic to $S_n$ if $X$ is a tree.
\end{lem}
\begin{proof}
Note that $\alpha(\epsilon)$ is the transposition $(E(x),E(y))$ for every $\epsilon=\{x,y\}\in E.$
Let $A=\{(E(x),E(y))\in S_n~|~\{x,y\}\in T\}.$ Pick any two distinct vertices $u,v\in V.$ Then there exists a subset $\{\{u_0,u_1\},\{u_1,u_2\},\ldots,\{u_{k-1},u_k\}\}$ of $T$ with $u_0=u$ and $u_k=v.$ Note that
\begin{align*}
(E(u),E(v))=&(E(u_{k-1}),E(u_{k}))\cdots(E(u_1),E(u_2))(E(u_0),E(u_1))\\
&(E(u_1),E(u_2))\cdots(E(u_{k-1}),E(u_{k})).
\end{align*}
Hence $A$ generates all transpositions in $S_n$ and then $A$ generates $S_n.$ Thus, the first assertion holds. For the second assertion, let
$X$ be a tree. Since $m=n-1,$ the edge space $\mathcal{E}$ is the bond space $\mathcal{B}$ of $X$ by Proposition~\ref{dimB}(ii). Hence, for $\mathbf{g}\in\mathbf{W}_E(X),$ if $\mathbf{g}E(v)=E(v)$ for every $v$ then $\mathbf{g}$ is the identity map on $\mathcal{E}.$ This shows that the kernel of $\alpha$ is trivial. From this and the first assertion, the second assertion holds.
\end{proof}



\begin{exam}
Let $X=(V,E)$ be the star graph of $n\geq 3$ vertices. By Lemma~\ref{WWTontoSn}, the edge-flipping group $\mathbf{W}_E(X)$ of $X$ is isomorphic to $S_n.$
\end{exam}

\section{Orbits}
Let $X=(V,E)$ denote a finite simple connected graph with $|V|=n$ and $|E|=m,$ and let $T$ denote a spanning tree of $E.$ Let $\mathbf{W}_E(X)$ denote the edge-flipping group of $X.$ In this section, we give a description of the orbits of the edge-flipping puzzle on $X$ in terms of our language. For this purpose, we fix a vertex $u$ in $V$ for this whole section and choose a nice basis of the bond space $\mathcal{B}$ of $X.$ Recall that from
Proposition~\ref{dimB}(iii) $E(u)=\sum_{v\in V-\{u\}}E(v),$ and from Proposition~\ref{dimB}(iv)
$$
\Delta:=\{E(v)~|~v\in V-\{u\}\}
$$
is a basis of $\mathcal{B}.$ We call $\Delta$ the {\it simple basis} of $\mathcal{B}.$ For each element $G$ in $\mathcal{B},$ let $\Delta(G)$ denote the subset of $\Delta$ such that the sum of all elements in $\Delta(G)$ is equal to $G,$ and let the {\it simple weight} $sw(G)$ of $G$ be the cardinality of $\Delta(G).$ For example, $\Delta(E(u))=\{E(v)~|~v\in V-\{u\}\}$ and
$sw(E(u))=n-1.$

Recall that an {\it orbit} of $\mathcal{E}$ under $\mathbf{W}_E(X)$ is the set $\{\mathbf{g}G~|~\mathbf{g}\in\mathbf{W}_E(X)\}$ for some $G\in\mathcal{E}.$ Note that the orbits of $\mathcal{E}$ under $\mathbf{W}_E(X)$ are corresponding to the orbits of the edge-flipping puzzle on $X$ and are the minimal
nonempty invariant subsets of $\mathcal{E}$ under $\mathbf{W}_E(X).$
Therefore, by Proposition~\ref{coset_rep} and Proposition~\ref{coset}, every orbit of $\mathcal{E}$ under $\mathbf{W}_E(X)$ is contained in $F+\mathcal{B}$ for some $F\subseteq E-T.$
In the following lemma, we give a description of the orbits of $\mathcal{B}$ under $\mathbf{W}_E(X)$ in terms of simple weights.

\begin{lem}(\cite[Theorem 10]{w:08})\label{orbit_under_WT}
Let $X=(V,E)$ be a finite simple connected graph with $|V|=n\geq 3$ and let $T$ be a spanning tree of $E.$ Then the orbits of $\mathcal{B}$ under $\mathbf{W}_E(X)$ are the same as the orbits of $\mathcal{B}$ under $\mathbf{W}_E(X)_T.$ More precisely, these orbits are
$$
\Omega_i:=\{G\in\mathcal{B}~|~sw(G)=i,n-i\} {\rm~for~}i=0,1,\ldots,\lceil \frac{n-1}{2}\rceil,
$$
where $\mathbf{W}_E(X)_T$ is the subgroup of $\mathbf{W}_E(X)$ generated by $\{\bm{\rho}_{\epsilon}~|~\epsilon\in T\}.$
\end{lem}
\noindent {\it Sketch of Proof.}
Recall that from Proposition~\ref{bond_1}(ii) and  Proposition~\ref{dimB}(i), the bond space $\mathcal{B}$ of $X$ consists of $E(U)=\sum_{v\in U}E(v)$ for all $U\subseteq V.$ By Lemma~\ref{WWTontoSn}, both $\mathbf{W}_E(X)$ and $\mathbf{W}_E(X)_T$ act on $\{E(v)~|~v\in V\}$ as the symmetric group on $\{E(v)~|~v\in V\}.$ Hence every orbit of $\mathcal{B}$ under $\mathbf{W}_E(X)$ (or $\mathbf{W}_E(X)_T$) is one of $\{E(U)~|~|U|=i\}$ for $0\leq i\leq n.$ Since $E(u)=\sum_{v\in V-\{u\}}E(v),$ both $\{E(U)~|~|U|=i\}$ and $\{E(U)~|~|U|=n-i\}$ are equal to $\Omega_i$ for $0\leq i\leq \lceil \frac{n-1}{2}\rceil.$         \hfill $\Box$

\bigskip

For nonempty $F\subseteq E-T,$ the orbits of $F+\mathcal{B}$ under $\mathbf{W}_E(X)$ in terms of simple weights is given in the following.

\begin{lem}(\cite[Theorem 12]{w:08})\label{orbit_of_coset_under_W}
Let $X=(V,E)$ be a finite simple connected graph with $|V|=n\geq 3$ and let $T$ be a spanning tree of $E.$ Let $F$ be a nonempty subset of $E-T$ and $\epsilon\in F.$ Then the orbits of $F+\mathcal{B}$ under $\mathbf{W}_E(X)$ are the same as the orbits of $F+\mathcal{B}$ under $\mathbf{W}_E(X)_{T\cup\{\epsilon\}}.$ More precisely, these orbits are
\begin{align}\label{orbit_of_coset_under_W_equation}
\left\{
\begin{array}{ll}
F+\mathcal{B}, &\hbox{if $n$ is odd;}\\
F+\mathcal{B}_e~{\rm and~} F+\mathcal{B}_o,  &\hbox{if $n$
is even,}
\end{array}
\right.
\end{align}
where $\mathbf{W}_E(X)_{T\cup\{\epsilon\}}$ is the subgroup of $\mathbf{W}_E(X)$ generated by $\{\bm{\rho}_{\epsilon'}~|~\epsilon'\in T\cup\{\epsilon\}\},$ $\mathcal{B}_e:=\{G\in \mathcal{B}~|~sw(G){\rm~is~even}\}$ and $\mathcal{B}_o:=\{G\in \mathcal{B}~|~sw(G){\rm~is~odd}\}.$
\end{lem}
\noindent {\it Sketch of Proof.}
We now determine the orbits of $F+\mathcal{B}$ under $\mathbf{W}_E(X)_{T\cup\{\epsilon\}}.$ Since $F\cap T=\emptyset$ and by (\ref{def_e}), $\bm{\rho}_{\epsilon'}F=F$ for any $\epsilon'\in T.$ Hence $\mathbf{W}_E(X)_TF=F,$ and the orbits of $F+\mathcal{B}$ under $\mathbf{W}_E(X)_T$ are
\begin{align}\label{orbit_of_coset_under_WT}
F+\Omega_i {\rm~for~} i=0,1,\ldots,\lceil \frac{n-1}{2}\rceil
\end{align}
by Lemma~\ref{orbit_under_WT}. It thus remains to consider the action of additional map $\bm{\rho}_{\epsilon}$ on $F+\mathcal{B}.$ Let $F+G\in F+\mathcal{B}$ with $G\in\mathcal{B}$ and $sw(G)=i$ for some $0\leq i\leq n-1.$
Note that $\bm{\rho}_{\epsilon}(F+G)=F+(E(\epsilon)+\bm{\rho}_{\epsilon}G)$ and  $E(\epsilon)+\bm{\rho}_{\epsilon}G\in\mathcal{B}.$ We discuss the simple weight of $E(\epsilon)+\bm{\rho}_{\epsilon}G$ as follows. If $u\notin \epsilon$ then $sw(E(\epsilon))=2$ and
\begin{align}\label{E(e)+eG1}
sw(E(\epsilon)+\bm{\rho}_{\epsilon}G)=\left\{
\begin{array}{ll}
i+2, &\hbox{if $|\Delta(G)\cap\Delta(E(\epsilon))|=0;$}\\
i, &\hbox{if $|\Delta(G)\cap\Delta(E(\epsilon))|=1;$}\\
i-2, &\hbox{otherwise,}
\end{array}
\right.
\end{align}
and if $u\in \epsilon$ then $sw(E(\epsilon))=n-2$ and
\begin{align}\label{E(e)+eG2}
sw(E(\epsilon)+\bm{\rho}_{\epsilon}G)=\left\{
\begin{array}{ll}
i, &\hbox{if $|\Delta(G)\cap\Delta(E(\epsilon))|=i-1;$}\\
n-i-2, &\hbox{otherwise.}\\
\end{array}
\right.
\end{align}
By (\ref{E(e)+eG1}) and (\ref{E(e)+eG2}), some sets in (\ref{orbit_of_coset_under_WT}) are further put together to become an orbit of $F+\mathcal{B}$ under $\mathbf{W}_E(X)_{T\cup\{\epsilon\}}$ and we have the result as described in (\ref{orbit_of_coset_under_W_equation}).
Since $\epsilon$ is an arbitrary element of $F,$ the orbits of $F+\mathcal{B}$ under $\mathbf{W}_E(X)_{T\cup F}$ and under $\mathbf{W}_E(X)_{T\cup \{\epsilon\}}$ are the same, where $\mathbf{W}_E(X)_{T\cup F}$ is the subgroup of $\mathbf{W}_E(X)$ generated by the set $\{\bm{\rho}_{\epsilon'}~|~\epsilon'\in T\cup F\}.$ Note that $\bm{\rho}_{\epsilon'}F=F$ for all $\epsilon'\in E-(T\cup F).$ Hence the orbits of $F+\mathcal{B}$ under $\mathbf{W}_E(X)$ and under $\mathbf{W}_E(X)_{T\cup F}$ are also the same.  \hfill $\Box$




\section{More on the edge-flipping groups}\label{Moreon_W}
Let $X=(V,E)$ denote a finite simple connected graph with $|V|=n\geq 3$ and $|E|=m.$ In this section, we investigate the structure of the edge-flipping group $\mathbf{W}_E(X)$ of $X.$ Let $\mathcal{B}_i$ be a copy of the bond space $\mathcal{B}$ of $X$ for $1\leq i\leq m-n+1.$ Recall that their {\it direct sum}
$$
\mathcal{B}^{m-n+1}:=\bigoplus^{m-n+1}_{i=1}\mathcal{B}_i
$$
is the set of all $(m-n+1)$-tuples $(G_i)^{m-n+1}_{i=1}$ where $G_i\in \mathcal{B}_i$ and where the addition is defined componentwise; i.e.,
$
(G_i)^{m-n+1}_{i=1}+(H_i)^{m-n+1}_{i=1}=(G_i+H_i)^{m-n+1}_{i=1}.
$
Let ${\rm Aut}(\mathcal{B}^{m-n+1})$ denote the automorphism group of $\mathcal{B}^{m-n+1}.$

\begin{defn}\label{beta}
Let $\beta:\mathbf{W}_E(X)\rightarrow {\rm Aut}(\mathcal{B}^{m-n+1})$ denote the group homomorphism defined by
$$
\beta(\mathbf{g})(G_i)^{m-n+1}_{i=1}=(\mathbf{g}G_i)^{m-n+1}_{i=1}
$$
for $\mathbf{g}\in \mathbf{W}_E(X)$ and $(G_i)^{m-n+1}_{i=1}\in\mathcal{B}^{m-n+1}.$
\end{defn}

Recall that from Lemma~\ref{WWTontoSn} the group homomorphism $\alpha$ from $\mathbf{W}_E(X)$ into the symmetric group $S_n$ on $\{E(v)~|~v\in V\}$ is surjective. The following lemma shows that there exists a unique group homomorphism $\theta:S_n\rightarrow {\rm Aut}(\mathcal{B}^{m-n+1})$ such that the diagram

\unitlength=1mm
\begin{picture}(120,40)
\put(45,33){\vector(1, 0){30}}  
\put(44,29){\vector(3, -2){31}}
\put(79,30){\vector(0, -1){19}}
\put(30,32){$\mathbf{W}_E(X)$}
\put(77,32){$S_n$}
\put(77,7){${\rm Aut}(\mathcal{B}^{m-n+1})$}
\put(60,34){$\alpha$}
\put(80,20){$\theta$}
\put(55,16){$\beta$}

\put(50,0){\footnotesize{{\bf Figure 1.}}}
\end{picture}

\noindent is commutative.

\begin{lem}\label{lem_semidirect}
There exists a unique group homomorphism $\theta:S_n\rightarrow {\rm Aut}$ $(\mathcal{B}^{m-n+1})$ such that $\beta=\theta\circ \alpha.$ Moreover, 
\begin{align}\label{theta}
\theta(\sigma)(E(v_{i}))^{m-n+1}_{i=1}=(\sigma(E(v_{i})))^{m-n+1}_{i=1}
\end{align}
for $v_1,v_2,\ldots,v_{m-n+1}\in V$ and $\sigma\in S_n.$
\end{lem}
\begin{proof}
Since $\alpha$ is surjective, if $\theta$ exists then $\theta$ is unique. To show the existence of $\theta,$ it suffices to show the kernel ${\rm Ker}~\alpha$ of $\alpha$ is contained in the kernel ${\rm Ker}~\beta$ of $\beta.$  If $\mathbf{g}\in{\rm Ker}~\alpha,$ then $\mathbf{g}E(v)=E(v)$ for all $v\in V$ and hence $\mathbf{g}\in {\rm Ker}~\beta$ since $\{E(v)~|~v\in V\}$ spans $\mathcal{B}.$ Pick $\sigma\in S_n$ and choose an element $\mathbf{h}$ in $\mathbf{W}_E(X)$ such that $\alpha(\mathbf{h})=\sigma.$ To prove~(\ref{theta}), it suffices to show that
$$
\beta(\mathbf{h})(E(v_{i}))^{m-n+1}_{i=1}=(\alpha(\mathbf{h})(E(v_{i})))^{m-n+1}_{i=1}
$$
for $v_1,v_2,\ldots,v_{m-n+1}\in V,$ since $\beta(\mathbf{h})=\theta(\sigma)$ and $\alpha(\mathbf{h})=\sigma.$ By Definition~\ref{alpha} and Definition~\ref{beta}, we obtain that both sides of the above equation are equal to $(\mathbf{h}E(v_1),\mathbf{h}E(v_2),$ $\ldots,\mathbf{h}E(v_{m-n+1}))$ as desired.
\end{proof}

In view of Lemma~\ref{lem_semidirect}, there is a {\it semidirect product} of $\mathcal{B}^{m-n+1}$ and $S_n$ with respect to $\theta$ \cite[p.155]{dj:02}, denoted by
$\mathcal{B}^{m-n+1}\rtimes_{\theta}S_n;$ i.e., $\mathcal{B}^{m-n+1}\rtimes_{\theta}S_n$ is the set $\mathcal{B}^{m-n+1}\times S_n$ with the group operation defined by
\begin{align}\label{semidirect}
\begin{split}
&((G_i)^{m-n+1}_{i=1},\sigma)~((H_i)^{m-n+1}_{i=1},\tau)\\
=&((G_i)^{m-n+1}_{i=1}+\theta(\sigma)(H_i)^{m-n+1}_{i=1},\sigma\tau)
\end{split}
\end{align}
for all $(G_i)^{m-n+1}_{i=1},(H_i)^{m-n+1}_{i=1}\in\mathcal{B}^{m-n+1}$ and $\sigma,\tau\in S_n.$ Recall that $T$ denotes a spanning tree of $E$ and $|T|=n-1.$ Let $E-T=\{\epsilon_1,\epsilon_2,\ldots,$ $\epsilon_{m-n+1}\}.$ Since $\{\epsilon_i\}+\mathbf{W}_E(X)\{\epsilon_i\}\subseteq \mathcal{B}$ for $1\leq i\leq m-n+1,$  we can define a map from $\mathbf{W}_E(X)$ into $\mathcal{B}^{m-n+1}\rtimes_{\theta}S_n$ as follows.

\begin{defn}\label{gamma}
Let $\gamma:\mathbf{W}_E(X)\rightarrow \mathcal{B}^{m-n+1}\rtimes_{\theta}S_n$ denote the map defined by
$$
\gamma(\mathbf{g})=((\{\epsilon_i\}+\mathbf{g}\{\epsilon_i\})^{m-n+1}_{i=1},\alpha(\mathbf{g}))
$$
for $\mathbf{g}\in\mathbf{W}_E(X).$
\end{defn}

The following lemma shows that $\gamma$ is a group monomorphism.

\begin{lem}\label{gamma_monomorphism}
$\gamma$ is a group monomorphism from $\mathbf{W}_E(X)$ into $\mathcal{B}^{m-n+1}\rtimes_{\theta}S_n.$
\end{lem}
\begin{proof}
For $\mathbf{g},\mathbf{h}\in \mathbf{W}_E(X),$
\begin{align*}
\gamma(\mathbf{g})\gamma(\mathbf{h})&=((\{\epsilon_i\}+\mathbf{g}\{\epsilon_i\})^{m-n+1}_{i=1},\alpha(\mathbf{g}))((\{\epsilon_i\}+\mathbf{h}\{\epsilon_i\})^{m-n+1}_{i=1},\alpha(\mathbf{h}))\\
&=((\{\epsilon_i\}+\mathbf{g}\{\epsilon_i\})^{m-n+1}_{i=1}+\theta(\alpha(\mathbf{g}))(\{\epsilon_i\}+\mathbf{h}\{\epsilon_i\})^{m-n+1}_{i=1},\alpha(\mathbf{g})\alpha(\mathbf{h}))\\
&=((\{\epsilon_i\}+\mathbf{g}\{\epsilon_i\})^{m-n+1}_{i=1}+\beta(\mathbf{g})(\{\epsilon_i\}+\mathbf{h}\{\epsilon_i\})^{m-n+1}_{i=1},\alpha(\mathbf{g}\mathbf{h}))\\
&=((\{\epsilon_i\}+\mathbf{g}\{\epsilon_i\})^{m-n+1}_{i=1}+(\mathbf{g}\{\epsilon_i\}+\mathbf{g}\mathbf{h}\{\epsilon_i\})^{m-n+1}_{i=1},\alpha(\mathbf{g}\mathbf{h}))\\
&=((\{\epsilon_i\}+\mathbf{g}\mathbf{h}\{\epsilon_i\})^{m-n+1}_{i=1},\alpha(\mathbf{g}\mathbf{h}))\\
&=\gamma(\mathbf{gh}).
\end{align*}
This shows that $\gamma$ is a group homomorphism. Since each $\mathbf{g}\in {\rm Ker}~\gamma$ fixes the spanning set $\{\{\epsilon_1\},\{\epsilon_2\},\ldots,\{\epsilon_{m-n+1}\}\}\cup\{E(v)~|~v\in V\}$ of the edge space $\mathcal{E}$ of $X,$ $\mathbf{g}$ is the identity map on $\mathcal{E}.$ Hence ${\rm Ker}~\gamma$ is trivial.
\end{proof}

By Lemma~\ref{gamma_monomorphism}, $\mathbf{W}_E(X)$ is isomorphic to the subgroup $\gamma(\mathbf{W}_E(X))$ of $\mathcal{B}^{m-n+1}\rtimes_{\theta}S_n.$ Fortunately, $\gamma(\mathbf{W}_E(X))$ is knowable. Recall that $\mathcal{B}_e,$ defined in Lemma~\ref{orbit_of_coset_under_W}, is an $(n-2)$-dimensional subspace of $\mathcal{B}.$ Let $\mathcal{B}_e^{m-n+1}$ denote the subgroup
$$
\bigoplus^{m-n+1}_{i=1}\mathcal{B}_{e,i}
$$
of $\mathcal{B}^{m-n+1},$
where $\mathcal{B}_{e,i}:=\mathcal{B}_e$ for $1\leq i\leq m-n+1.$

\begin{thm}\label{W_iso_thm}
Let $X=(V,E)$ be a finite simple connected graph with $n\geq 3$ vertices and $m$ edges. Then the edge-flipping group $\mathbf{W}_E(X)$ is isomorphic to
\begin{align*}
\left\{
\begin{array}{ll}
\mathcal{B}^{m-n+1}\rtimes_{\theta}S_n, &\hbox{if $n$ is odd;}\\
\mathcal{B}^{m-n+1}_e\rtimes_{\theta}S_n, &\hbox{if $n$ is even.}
\end{array}
\right.
\end{align*}
\end{thm}
\begin{proof}
It suffices to show that for any $\sigma\in S_n,$ there exists $\mathbf{g}\in \mathbf{W}_E(X)$ such that
\begin{align}\label{gamma(W_T)}
\gamma(\mathbf{g})=((\emptyset)^{m-n+1}_{i=1},\sigma),
\end{align}
and for each $1\leq i\leq m-n+1,$ for any
$$
G\in\left\{
\begin{array}{ll}
B_i, &\hbox{if $n$ is odd;}\\
B_{e,i} &\hbox{if $n$ is even,}
\end{array}
\right.
$$
there exists $\mathbf{h}\in \mathbf{W}_E(X)$ such that
\begin{align}\label{gamma(W_Te)}
\gamma(\mathbf{h})=(\emptyset,\ldots,\emptyset, G,\emptyset,\ldots,\emptyset,\alpha(\mathbf{h})),
\end{align}
where $G$ is in the $i$th coordinate. From Lemma \ref{WWTontoSn}, (\ref{gamma(W_T)}) follows by choosing $\mathbf{g}\in \mathbf{W}_E(X)_T$ with $\alpha(\mathbf{g})=\sigma,$ since $\mathbf{g}\{\epsilon_j\}=\{\epsilon_j\}$ for each $1\leq j\leq m-n+1.$ From Lemma \ref{orbit_of_coset_under_W}, (\ref{gamma(W_Te)}) follows by choosing $\mathbf{h}\in\mathbf{W}_E(X)_{T\cup\{\epsilon_i\}}$ with $\mathbf{h}\{\epsilon_i\}=\{\epsilon_i\}+G,$ since  $\mathbf{h}\{\epsilon_j\}=\{\epsilon_j\}$ for $j\not= i.$
\end{proof}

Since ${\rm dim}~\mathcal{B}=n-1$ and ${\rm dim}~\mathcal{B}_e=n-2,$ the additive groups of $\mathcal{B}$ and $\mathcal{B}_e$ are isomorphic to $(\mathbb{Z}/2\mathbb{Z})^{n-1}$ and $(\mathbb{Z}/2\mathbb{Z})^{n-2}$ respectively, where $\mathbb{Z}$ is the additive group of integers.

\begin{exam}
Let $X$ be a cycle of $n$ vertices. Then the edge-flipping group $\mathbf{W}_E(X)$ of $X$ is isomorphic to
$$
\left\{
\begin{array}{ll}
(\mathbb{Z}/2\mathbb{Z})^{n-1}\rtimes S_n, & \hbox{if $n$ is odd;}\\
(\mathbb{Z}/2\mathbb{Z})^{n-2}\rtimes S_n, & \hbox{if $n$ is even}
\end{array}
\right.
$$
by Theorem \ref{W_iso_thm}.
\end{exam}

The following corollary says that there is a unique (up to isomorphism) edge-flipping group $\mathbf{W}_E(X)$ of all finite simple connected graphs $X=(V,E)$ with $|V|=n\geq 3$ and $|E|=m.$

\begin{cor}\label{unique_W_linegraphs}
Let $X=(V,E)$ and $X'=(V',E')$ be two finite simple connected graphs with $|V|=|V'|=n\geq 3$ and $|E|=|E'|.$  Then the edge-flipping group $\mathbf{W}_{E}(X)$ of $X$ is isomorphic to the edge-flipping group $\mathbf{W}_{E'}(X')$ of $X'.$
\end{cor}
\begin{proof}
We may assume that $V'=V.$ Recall that the simple basis $\Delta$ of $\mathcal{B}$ is the set $\{E(v)~|~v\in V-\{u\}\}$ for some fixed vertex $u\in V.$  Define $E'(v),~\mathcal{B}',~\Delta':=\{E'(v)~|~v\in V-\{u\}\},~\mathcal{B}_e',~S_n'$ and $\theta'$ correspondingly. From Theorem \ref{W_iso_thm}, it suffices to show that $\mathcal{B}^{m-n+1}\rtimes_{\theta}S_n$ and $\mathcal{B}_e^{m-n+1}\rtimes_{\theta}S_n$ are isomorphic to $\mathcal{B'}^{m-n+1}\rtimes_{\theta'} S_n'$ and $\mathcal{B}_e'^{m-n+1}\rtimes_{\theta'} S_n'$ respectively. Let  $\mu:\mathcal{B}\rightarrow\mathcal{B}'$ denote the invertible linear transformation defined by
$$
\mu(E(v))=E'(v)
$$
for $v\in V-\{u\}.$ Note that there exists a unique group isomorphism $\mu_{*}:S_n\rightarrow S_n'$ such that
$$
\mu_{*}(\sigma)(E'(v))=\mu(\sigma(E(v)))
$$
for all $\sigma\in S_n$ and $v\in V.$ Hence, there exists a unique bijective map $\phi:\mathcal{B}^{m-n+1}\rtimes_{\theta}S_n\rightarrow \mathcal{B'}^{m-n+1}\rtimes_{\theta'} S_n'$ such that
$$
\phi((G_i)^{m-n+1}_{i=1},\sigma)=((\mu(G_i))^{m-n+1}_{i=1},\mu_{*}(\sigma))
$$
for all $(G_i)^{m-n+1}_{i=1}\in \mathcal{B}^{m-n+1}$ and $\sigma\in S_n.$ The map $\phi$ sends $\mathcal{B}_e^{m-n+1}\rtimes_{\theta} S_n$ to $\mathcal{B}_e'^{m-n+1}\rtimes_{\theta'} S_n'$ since $\mu(\mathcal{B}_e)=\mathcal{B}_e'.$ By (\ref{theta}) and (\ref{semidirect}), it is straightforward to verify that $\phi$ is a group isomorphism, as desired.
\end{proof}

\section{Applications}
Let $X=(V,E)$ denote a finite simple connected graph with $|E|=m.$ In this section, we investigate the vertex-flipping group $\mathbf{W}_V(X)$ of $X,$ which is a vertex version of the edge-flipping group $\mathbf{W}_E(X)$ of $X.$ The vertex-flipping group $\mathbf{W}_V(X)$ of $X$ is also a subgroup of the general linear group ${\rm GL}(\mathcal{V})$ of the vertex space $\mathcal{V}$ of $X.$ The orbits of $\mathcal{V}$ under $\mathbf{W}_V(X)$ are corresponding to the orbits of the vertex-flipping puzzle on $X.$ See \cite{hw:pre,xw:07} for details. For $v\in V,$ let $N(v)$ denote the set consisting of all neighbors of $v.$ In the following, we give a formal definition of $\mathbf{W}_V(X).$


\begin{defn}\label{vertex-flipping_gp}
The {\it vertex-flipping group} $\mathbf{W}_V(X)$ of $X=(V,E)$ is the subgroup of the general linear group ${\rm GL}(\mathcal{V})$ of $\mathcal{V}$ generated by $\{\bm{s}_v~|~v\in V\},$ where $\bm{s}_v:\mathcal{V}\rightarrow\mathcal{V}$ is defined by
\begin{align*}
\bm{s}_vU=\left\{
\begin{array}{ll}
U+N(v), &\hbox{if $v\in U;$}\\
U,  &\hbox{otherwise}
\end{array}
\right.
\end{align*}
for $U\in \mathcal{V}.$ 
\end{defn}

The {\it line graph} $L(X)$ of $X=(V,E)$ is a finite simple connected graph with vertex set $E$ and edge set  $\{\{\epsilon,\epsilon'\}~|~|\epsilon\cap\epsilon'|=1{\rm~for~}\epsilon,\epsilon'\in E\}.$ From this definition, the edge-flipping group of $X$ and the vertex-flipping group of $L(X)$ are the same.


For this whole section, let $Y=(Z,H)$ denote the finite simple connected graph with $Z=\{0,1,2,\ldots$ $,m-1\}$ and $H=\{\{1,2\},\{2,3\},\ldots,\{m-2,m-1\},\{0,i_1\},\{0,i_2\},\ldots,\{0,i_\ell\}\},$ where $m\geq 2$ and $1\leq i_1<i_2<\cdots<i_\ell\leq m-1.$ See Figure 2. For example, the simply-laced Dynkin diagrams and extended Dynkin diagrams $\widetilde{A}_n,\widetilde{E}_7,\widetilde{E}_8$ are such graphs.


\unitlength=1mm
\begin{picture}(120,35)
\multiput(30,10)(6, 0){2}{\line(1, 0){6}} 
\multiput(90,10)(-6,0){2}{\line(-1, 0){6}}
\qbezier[180](48,10)(54,20)(60,30)
\qbezier[180](72,10)(66,20)(60,30)
\qbezier[180](66,10)(63,20)(60,30)
\multiput(30,10)(6, 0){4}{\circle*{1}} 
\multiput(90,10)(-6,0){4}{\circle*{1}}
\multiput(66,10)(-6,0){1}{\circle*{1}}
\multiput(60,30)(6,0){1}{\circle*{1}}
\multiput(42,10)(3,0){12}{\circle*{0.5}} 
\multiput(61,20)(-2,0){4}{\circle*{0.5}}
\put(27,8){\tiny{$m-1$}} 
\put(34,8){\tiny{$m-2$}}
\put(49,8){\tiny{$i_\ell$}}
\put(91,8){\tiny{$1$}}
\put(85,8){\tiny{$2$}}
\put(79,8){\tiny{$3$}}
\put(73,8){\tiny{$i_1$}}
\put(67,8){\tiny{$i_2$}}
\put(61,30){\tiny{$0$}}
\put(43,0){\footnotesize{{\bf Figure 2.} The graph $Y$}}
\end{picture}
\bigskip

Let the value $\pi_1$ of $Y$ be
\begin{align}\label{value_pi1}
\left\{
\begin{array}{ll}
\sum\limits^\ell_{t=1}(-1)^ti_{t}, & \hbox{if $\ell$ is even;}\\
\sum\limits^\ell_{t=1}(-1)^ti_{t}+m, & \hbox{otherwise.}
\end{array}
\right.
\end{align}
Note that $1\leq \pi_1\leq m-1.$

\begin{thm}(\cite[Theorem 3.9]{hw:pre2})\label{Theorem 3.9} Let $1\leq k\leq m-1$ be an integer. Then the vertex-flipping group of $Y=(Z,H)$ is unique (up to isomorphism) among those graphs $Y$ with $\pi_1=k.$ \hfill $\Box$
\end{thm}

The aim of this section is to determine the structures of the vertex-flipping groups $\mathbf{W}_Z(Y)$ of some graphs $Y.$ For this purpose, we define some terms. Two finite simple graphs $X=(V,E)$ and $X'=(V',E')$ are {\it isomorphic} if there exists a bijective map $\phi:V\rightarrow V'$ such that $\{x,y\}\in E$ if and only if $\{\phi(x),\phi(y)\}\in E'$ for all $x,y\in V.$ We shall denote that two finite simple graphs $X$ and $X'$ are isomorphic by writing $X\cong X'.$ 
A graph $(U,F)$ is a {\it subgraph} of $X=$ $(V,E)$ if $U\subseteq V$ and $F\subseteq E,$ and a subgraph $(U,F)$ of $X=(V,E)$ is {\it induced} if $F=\{\{x,y\}\in E~|~x,y\in U\}.$ An (induced) subgraph $(U,F)$ of $X=(V,E)$ is an {\it (induced) path} if $U=\{v_0,v_1,\ldots,v_k\}$ and $F=\{\{v_0,v_1\},\{v_1,v_2\},\ldots,\{v_{k-1},v_k\}\}$ for some nonnegative integer $k,$ where $v_i$ are all distinct. The following easy fact will be used later.

\begin{lem}\label{induced_path}
Let $X=(V,E)$ be a finite simple connected graph with $|E|=m.$ 
Let $1\leq k\leq m$ be an integer.
Then $X$ contains a path of $k$ edges if and only if the line graph $L(X)$ of $X$ contains an induced path of $k$ vertices. \hfill $\Box$
\end{lem}

We determine the structure of the vertex-flipping group $\mathbf{W}_Z(Y)$ of $Y=(Z,H)$ in some special cases. 

\begin{cor}\label{W_pi_1-1,2,n-2,n-1}
Let $Y=(Z,H)$ be a finite simple connected graph with $Z=\{0,1,2,\ldots$ $,m-1\}$ and $H=\{\{1,2\},\{2,3\},\ldots,\{m-2,m-1\},\{0,i_1\},\{0,i_2\},$ $\ldots,\{0,i_\ell\}\},$ where $1\leq i_1<i_2<\cdots<i_\ell\leq m-1.$ Let the value $\pi_1$ of $Y$ be defined in~(\ref{value_pi1}). Then the following (i)-(iii) hold.
\begin{enumerate}
\item[(i)] If $Y$ is isomorphic to a line graph $L(X)$ for some finite simple connected graph $X,$ then $\pi_1\in\{1,2,m-2,m-1\}.$

\item[(ii)] If $\pi_1\in\{1,m-1\},$ then the vertex-flipping group $\mathbf{W}_Z(Y)$ of $Y$ is isomorphic to the symmetric group $S_{m+1}$ of degree $m+1.$

\item[(iii)] If $\pi_1\in\{2,m-2\},$ then
$\mathbf{W}_Z(Y)$ is isomorphic to
$$
\left\{
\begin{array}{ll}
(\mathbb{Z}/2\mathbb{Z})^{m-1}\rtimes S_m, &\hbox{if $m$ is odd;}\\
(\mathbb{Z}/2\mathbb{Z})^{m-2}\rtimes S_m, &\hbox{if $m$ is even,}
\end{array}
\right.
$$
where $\mathbb{Z}$ is the additive group of integers.
\end{enumerate}
\end{cor}
\begin{proof}
Suppose that $Y$ is isomorphic to $L(X)$ for some finite simple connected graph $X.$ Since $Y$ contains the induced path
$$
(\{1,2,\ldots,m-1\},\{\{1,2\},\{2,3\},\ldots,\{m-2,m-1\}\}),
$$
$X$ contains a path of $m-1$ edges by Lemma~\ref{induced_path}. Note that $X$ has $m$ edges, one more edge besides the path. Hence the left column of Figure 3 completely lists all such graphs $X,$ and the right column is the corresponding line graph $L(X)$ of $X.$ By computing the value $\pi_1$ of $Y\cong L(X)$ and using (\ref{value_pi1}), we find $\pi_1=1,2,m-2,$ or $m-1.$ This proves (i). Note that the vertex-flipping group $\mathbf{W}_Z(Y)$ of $Y$ is the edge-flipping group of $X,$ and its group structure only depends on $\pi_1$ by Theorem~\ref{Theorem 3.9}. Hence we find $\mathbf{W}_Z(Y)$ as listed in (ii),(iii) by Theorem~\ref{W_iso_thm}.

\begin{picture}(125, 129)
\put(28,121){$X$}
\put(58,123){\vector(1, 0){9}}
\put(85,121){$Y\cong L(X)$}
\multiput(10,104)(5, 0){3}{\line(1, 0){5}}
\multiput(50,104)(-5,0){3}{\line(-1, 0){5}}
\multiput(10,104)(0,5){1}{\line(0,1){5}}
\multiput(10,104)(5, 0){4}{\circle*{1}} 
\multiput(50,104)(-5,0){4}{\circle*{1}}
\multiput(10,109)(-5,0){1}{\circle*{1}}
\multiput(27,104)(2,0){4}{\circle*{0.5}} 
\put(8,106){\scriptsize $0$} 

\put(58, 104){\vector(1, 0){9}}

\multiput(75, 104)(5, 0){3}{\line(1, 0){5}} 
\multiput(115,104)(-5,0){3}{\line(-1, 0){5}}
\qbezier[250](75,104)(85,109)(95,114)
\multiput(75, 104)(5, 0){4}{\circle*{1}} 
\multiput(115,104)(-5,0){4}{\circle*{1}}
\multiput(95,114)(-5,0){1}{\circle*{1}}
\multiput(92,104)(3,0){3}{\circle*{0.5}} 
\put(96,114){\scriptsize $0$} 

\put(84,98){\footnotesize{$(\pi_1=1,m-1)$}}

\multiput(50,82)(-5,0){2}{\line(-1, 0){5}}
\multiput(18,82)(5,0){2}{\line(1, 0){5}}
\multiput(23,82)(0,5){1}{\line(0,1){5}}
\multiput(10,82)(5,0){1}{\circle*{1}} 
\multiput(18,82)(5,0){3}{\circle*{1}}
\multiput(50,82)(-5,0){3}{\circle*{1}}
\multiput(23,87)(0,5){1}{\circle*{1}}
\multiput(12,82)(2,0){3}{\circle*{0.5}} 
\multiput(30,82)(2,0){5}{\circle*{0.5}}
\put(24,84){\scriptsize $0$} 

\put(58, 82){\vector(1, 0){9}}

\multiput(115,82)(-5, 0){2}{\line(-1, 0){5}} 
\multiput(83,82)(5,0){1}{\line(1,0){5}}
\qbezier[100](83,82)(89,87)(95,92)
\qbezier[100](88,82)(91.5,87)(95,92)
\multiput(115,82)(-5, 0){3}{\circle*{1}} 
\multiput(83,82)(5,0){2}{\circle*{1}}
\multiput(75,82)(5,0){1}{\circle*{1}}
\multiput(95,92)(5,0){1}{\circle*{1}}
\multiput(102.5,82)(-3,0){6}{\circle*{0.5}} 
\multiput(77.5,82)(3,0){2}{\circle*{0.5}}
\put(96,92){\scriptsize $0$} 

\put(89,76){\footnotesize{$(\pi_1=1)$}}

\qbezier[500](10,60)(30,80)(50,60)
\multiput(10,60)(5, 0){3}{\line(1, 0){5}}
\multiput(50,60)(-5,0){3}{\line(-1, 0){5}}
\multiput(10,60)(5, 0){4}{\circle*{1}} 
\multiput(50,60)(-5,0){4}{\circle*{1}}
\multiput(27,60)(2,0){4}{\circle*{0.5}} 
\put(30,71){\scriptsize $0$} 

\put(58, 60){\vector(1, 0){9}}

\multiput(75, 60)(5, 0){3}{\line(1, 0){5}} 
\multiput(115,60)(-5,0){3}{\line(-1, 0){5}}
\qbezier[250](75,60)(85,65)(95,70)
\qbezier[250](115,60)(105,65)(95,70)
\multiput(75, 60)(5, 0){4}{\circle*{1}} 
\multiput(115,60)(-5,0){4}{\circle*{1}}
\multiput(95,70)(-5,0){1}{\circle*{1}}
\multiput(92,60)(3,0){3}{\circle*{0.5}} 
\put(96,70){\scriptsize $0$} 

\put(85.5,54){\footnotesize{$(\pi_1=m-2)$}}

\qbezier[400](10,38)(23.5,58)(37,38)
\multiput(10,38)(5,0){2}{\line(1, 0){5}}
\multiput(42,38)(-5,0){2}{\line(-1, 0){5}}
\multiput(10,38)(5,0){3}{\circle*{1}} 
\multiput(42,38)(-5,0){3}{\circle*{1}}
\multiput(50,38)(-5,0){1}{\circle*{1}}
\multiput(22,38)(2,0){5}{\circle*{0.5}} 
\multiput(48,38)(-2,0){3}{\circle*{0.5}}
\put(23,49){\scriptsize $0$} 

\put(58, 38){\vector(1, 0){9}}

\multiput(75, 38)(5, 0){2}{\line(1, 0){5}} 
\multiput(107,38)(-5,0){1}{\line(-1,0){5}}
\qbezier[250](75,38)(85,43)(95,48)
\qbezier[100](107,38)(101,43)(95,48)
\qbezier[100](102,38)(98.5,43)(95,48)
\multiput(75, 38)(5, 0){3}{\circle*{1}} 
\multiput(107,38)(-5,0){2}{\circle*{1}}
\multiput(115,38)(-5,0){1}{\circle*{1}}
\multiput(95,48)(-5,0){1}{\circle*{1}}
\multiput(87.5,38)(3,0){6}{\circle*{0.5}} 
\multiput(112.5,38)(-3,0){2}{\circle*{0.5}}
\put(96,48){\scriptsize $0$} 

\put(84,32){\footnotesize{$(\pi_1=2,m-2)$}}

\qbezier[250](22,16)(30,36)(38,16)
\multiput(18,16)(4, 0){2}{\line(1, 0){4}}
\multiput(42,16)(-4, 0){2}{\line(-1, 0){4}}
\multiput(10,16)(4, 0){1}{\circle*{1}} 
\multiput(18,16)(4, 0){3}{\circle*{1}}
\multiput(42,16)(-4, 0){3}{\circle*{1}}
\multiput(50,16)(-4,0){1}{\circle*{1}}
\multiput(12,16)(2,0){3}{\circle*{0.5}} 
\multiput(28,16)(2,0){3}{\circle*{0.5}}
\multiput(48,16)(-2,0){3}{\circle*{0.5}}
\put(30,27){\scriptsize $0$} 

\put(58, 16){\vector(1, 0){9}}

\multiput(83, 16)(5, 0){1}{\line(1, 0){5}} 
\multiput(107,16)(-5,0){1}{\line(-1,0){5}}
\qbezier[100](83,16)(89,21)(95,26)
\qbezier[100](88,16)(91.5,21)(95,26)
\qbezier[100](107,16)(101,21)(95,26)
\qbezier[100](102,16)(98.5,21)(95,26)
\multiput(75, 16)(5, 0){1}{\circle*{1}} 
\multiput(83, 16)(5, 0){2}{\circle*{1}}
\multiput(107,16)(-5,0){2}{\circle*{1}}
\multiput(115,16)(-5,0){1}{\circle*{1}}
\multiput(95,26)(-5,0){1}{\circle*{1}}
\multiput(77.5,16)(3,0){2}{\circle*{0.5}} 
\multiput(90.5,16)(3,0){4}{\circle*{0.5}}
\multiput(112.5,16)(-3,0){2}{\circle*{0.5}}
\put(96,26){\scriptsize $0$} 

\put(89,10){\footnotesize{$(\pi_1=2)$}}
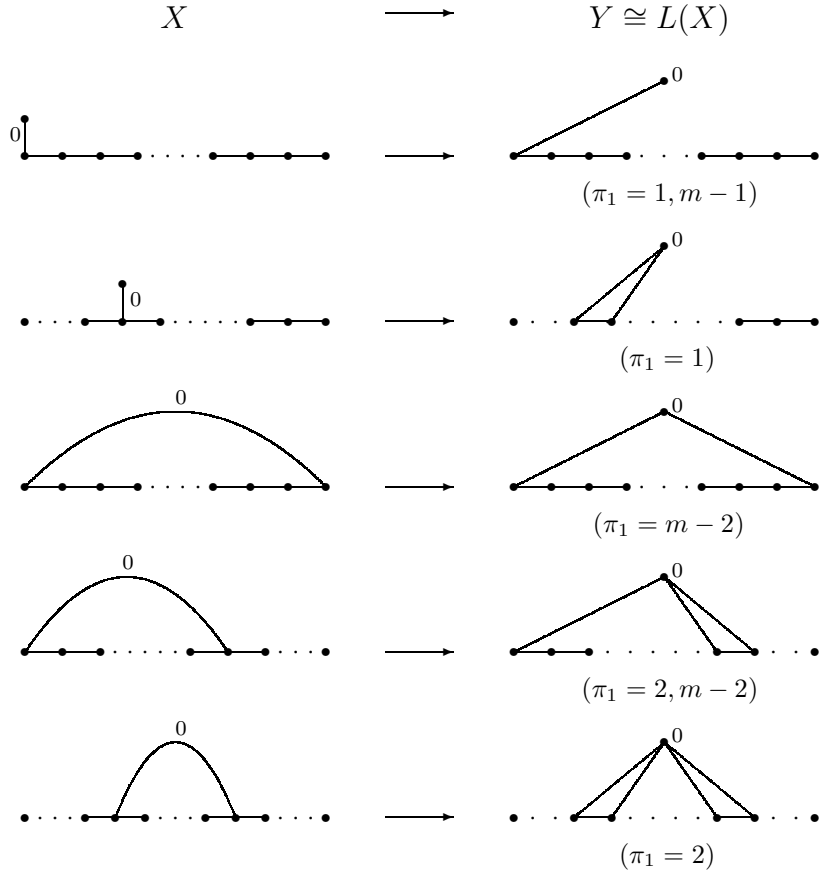
\put(22,3){\footnotesize{{\bf Figure 3.} All graphs $Y$ are isomorphic to line graphs}}
\end{picture}
\end{proof}

\begin{exam}
The graph $Y=(Z,H)$ in Figure 4

\begin{picture}(120,31)
\multiput(48,12)(8,0){3}{\line(1,0){8}} 
\qbezier[150](48,12)(54,19.5)(60,27)
\qbezier[150](64,12)(62,19.5)(60,27)
\multiput(48,12)(8,0){4}{\circle*{1}} 
\put(47,9){\scriptsize $4$}
\put(55,9){\scriptsize $3$}
\put(63,9){\scriptsize $2$}
\put(71,9){\scriptsize $1$}
\multiput(60,27)(8,0){1}{\circle*{1}}
\put(58,27){\scriptsize $0$}
\put(54,2){\footnotesize{\bf{Figure 4.}}}
\end{picture}

\noindent is a five-vertex graph containing an induced path of four vertices. By (\ref{value_pi1}), its value $\pi_1$ is 2. Hence, $\mathbf{W}_Z(Y)$ is isomorphic to $(\mathbb{Z}/2\mathbb{Z})^4\rtimes S_5$ by Corollary \ref{W_pi_1-1,2,n-2,n-1}(iii).
\end{exam}

In the case of $\pi_1\in\{1,2,m-2,m-1\},$ we use (\ref{value_pi1}) to find all such graphs, and there are only two graphs that are not isomorphic to line graphs. We show both of them in Figure 5. Note that their possible values $\pi_1$ are in $\{2,m-2\}.$ By Corollary \ref{W_pi_1-1,2,n-2,n-1}(iii), if $Y=(Z,H)$ is one of two graphs in Figure 5, the vertex-flipping group $\mathbf{W}_Z(Y)$ of $Y$ is isomorphic to
\begin{align*}\label{W_notlinegraph}
\left\{
\begin{array}{ll}
(\mathbb{Z}/2\mathbb{Z})^{m-1}\rtimes S_m, &\hbox{if $m$ is odd;}\\
(\mathbb{Z}/2\mathbb{Z})^{m-2}\rtimes S_m, &\hbox{if $m$ is even.}
\end{array}
\right.
\end{align*}

\begin{picture}(120,30)

\multiput(10,14)(5,0){3}{\line(1,0){5}} 
\multiput(50,14)(-5,0){3}{\line(-1,0){5}}
\qbezier[180](15,14)(22.5,19)(30,24)
\multiput(10,14)(5,0){4}{\circle*{1}} 
\multiput(50,14)(-5,0){4}{\circle*{1}}
\multiput(30,24)(-6,0){1}{\circle*{1}}
\multiput(27,14)(3,0){3}{\circle*{0.5}} 

\multiput(90,14)(5,0){3}{\line(1,0){5}} 
\multiput(120,14)(-5,0){1}{\line(-1,0){5}}
\qbezier[120](90,14)(95,19)(100,24)
\qbezier[80](100,14)(100,19)(100,24)
\multiput(80,14)(5,0){1}{\circle*{1}} 
\multiput(90,14)(5,0){4}{\circle*{1}}
\multiput(120,14)(-5,0){2}{\circle*{1}}
\multiput(100,24)(-6,0){1}{\circle*{1}}
\multiput(82,14)(3,0){3}{\circle*{0.5}} 
\multiput(113,14)(-3,0){3}{\circle*{0.5}}

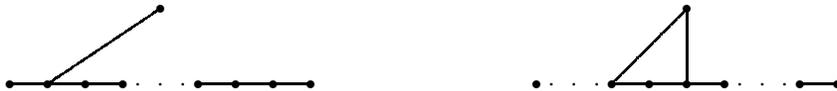
\put(25,5){\footnotesize{{\bf Figure 5.} All graphs $Y$ with $\pi_1\in\{1,2,m-2,m-1\}$}}
\put(43,1){\footnotesize{are not isomorphic to line graphs}}
\end{picture}

\begin{rem}\label{r5.8}
Theorem~\ref{Theorem 3.9} implies that in the class of $2^{m-1}$ graphs $Y=(Z,H),$ the number of the vertex-flipping groups $\mathbf{W}_Z(Y)$ of $Y$ is at most $m-1$ up to isomorphism. Together with Corollary~\ref{unique_W_linegraphs}, it seems that for a given vertex number, the number of non-isomorphic vertex-flipping groups is not too large, and a classification of them seems to be visible.
\end{rem}

\section{Acknowledgements}
The authors thank the anonymous referees for giving many valuable suggestions
in the presentation of the paper.

\bigskip

\noindent Hau-wen Huang \hfil\break Department of Applied Mathematics \hfil\break National Chiao
Tung University \hfil\break 1001 Ta Hsueh Road \hfil\break Hsinchu, Taiwan 30050, R.O.C.\hfil\break
Email: {\tt poker80.am94g@nctu.edu.tw} \hfil\break Fax: +886-3-5724679 \hfil\break
\medskip

\noindent Chih-wen Weng \hfil\break Department of Applied Mathematics \hfil\break National Chiao
Tung University \hfil\break 1001 Ta Hsueh Road \hfil\break Hsinchu, Taiwan 30050, R.O.C.\hfil\break
Email: {\tt weng@math.nctu.edu.tw} \hfil\break Fax: +886-3-5724679 \hfil\break
\medskip

\end{document}